\documentclass[11pt,a4paper]{article}

\usepackage{epic,eepic,epsfig,amsmath,amssymb,amsthm}
\usepackage{t1enc}

\usepackage{color}

\usepackage{bbm}

\def\cF{{\mathcal F}}

\def\virgp{\raise 2pt\hbox{,}}

\renewcommand{\geq}{\geqslant}
\renewcommand{\leq}{\leqslant}
\def\N{{\mathbb N}}

\def\R{{\mathbb R}}

\def\virgp{\raise 2pt\hbox{,}}
\def\cdotpv{\raise 2pt\hbox{;}}

\def\1{\mathbbm{1}}

\newtheorem{theorem}{Theorem}[section]

\newtheorem{proposition}[theorem]{Proposition}

\theoremstyle{remark}
\newtheorem{remark}{Remark}[section]

\theoremstyle{definition}

\theoremstyle{definition}

\theoremstyle{definition}

\begin{document}

\title{Control of the Black-Scholes equation}

\author{Claire David}

\maketitle
\centerline{Universit\'e Pierre et Marie Curie-Paris 6}
\centerline{Laboratoire Jacques Louis Lions - UMR 7598}
\centerline{Boîte courrier 187, 4 place Jussieu, F-75252 Paris
cedex 05, France}

\begin{abstract}

The purpose of this work is to apply the results developed by J.Y. Chemin and Cl. David \cite{CheminEtClaire1}, \cite{CheminEtClaire2}, to the Black-Scholes equation. This latter equation being directly linked to the heat equation, it enables us to propose a new approach allowing to control properties of the solution by means of a shape parameter.
\end{abstract}



\section{Introduction}

\noindent It is well known that the Black-Scholes ($\cal {BS}$) model, which gives the dynamic of the option prices in financial markets \cite{BlackScholes}, \cite{Merton}, is given, in its usual arbitrage free version, and the specific case of a "call", for european options, by:

$$\partial_t C + \displaystyle \frac{\sigma^2 \,S^2}{2} \displaystyle \frac{\partial^2 C}{\partial S^2}+ r\,\left ( S \,\partial_S C - C \right)= 0$$

\noindent where the price of the option, $C$, is a function of the underlying asset price $S$ and time $t$; $ r$ is the risk-free interest rate,
and $ \sigma$ the volatility of the stock.\\

\noindent There are hundreds of papers dealing with the Black-Scholes equation. Yet, no one seems to have ever used the scaling invariance coming from the heat equation. And there appears to be very few studies on the "control aspect" of the equation (one can see, for instance \cite{FreyPolte}).\\

\noindent  In \cite{CheminEtClaire2}, J.Y. Chemin and Cl. David obtained new results for the mass critical nonlinear Schr\"{o}dinger equation, building a continuous map from the Lebesgue space~$L^2(\R^2)$ in the set~$\cal G$ of initial data which give birth to global solution in the space~$L^4(\R^{1+2})$. The principle is simple: one uses the fact that for this nonlinear equation, solutions of scales which are different enough almost do not interact. The nonlinearity of the original equation just required to determine a condition about the size of the scale which depends continuously on the data. \\

\noindent The idea is to apply the building technique of the map described above, starting from the point that the Black-Scholes equation can be transformed into a linear heat equation, which has a natural scaling invariance, and that the linearity of the transformed equation enables one to obtain exact solutions, and, thus, yield interesting results. The aim of this paper is to present a new way to change either the Delta, or other greeks, by a given amount ; this technique appears thus as a new way to control the equation, without resorting to control functions. Indeed, one can change the greeks directly using the general solution of
the Black-Scholes equation, nevertheless, our technique appears as an alternative and simple one to implement.

\noindent Section~\ref{Build} is devoted to applying the Chemin-David technique to the Black-Scholes equation.\\

\noindent The consistency of the approach with respect to initial conditions is examined in Section~\ref{Consistency}.\\

\noindent Section~\ref{Results} present numerical results.\\

\section{Introduction of a shape parameter}

\label{Build}

\subsection{Classical results for the Black-Scholes equation}

\noindent If one denotes by $E$ the exercise price of the option, for a study in the time interval $[0,T]$, $T > 0$, the limit conditions are:

$$\left \lbrace \begin{array}{ccc} C(0, t)  &=& 0\quad \text{ for all } t \\
 \displaystyle \lim_{S \to + \infty} C(S, t) &=& S\\
   C(S, T)  &=& \max\, \left \lbrace S - E, 0 \right \rbrace \end{array} \right.$$

\noindent The Black-Scholes equation can classically be transformed into a heat-diffusion equation, through the change of variables:


$$\tau=   \displaystyle\frac{\sigma^2}{2 } \,\left (T-t\right)  \quad , \quad
  x = \ln \displaystyle\frac{S}{E}
 $$

\noindent setting:

$$  C(S,t)=E\, \widetilde{C}(x,\tau)$$

\noindent and:

 $$ \widetilde{C} (x,\tau)= e^{ \alpha\,x+\beta\,\tau}\, \widetilde{\widetilde{C}}(x,\tau)$$

\noindent where:

$$k=  \displaystyle\frac{2\,r }{\sigma^2} \quad , \quad \alpha= \displaystyle  \frac{1-k}{2} \quad , \quad
 \beta=\alpha^2+(k-1)\,\alpha-k=-
  \displaystyle  \frac{(k+1)^2}{4}
  $$

\noindent which lead to the normalized heat equation:

\begin{equation}
\label{Heat}\displaystyle  \frac{\partial \widetilde{\widetilde{C}}}{\partial\tau} = \displaystyle  \frac{\partial^2 \widetilde{\widetilde{C}}}{\partial x^2}
\qquad \forall \,(x,\tau)\,\in\,\left [0,\displaystyle  \frac{\sigma^2\,T}{2}\right] \times \R
\end{equation}

\noindent with the initial condition:

\begin{equation}\label{CondIn2} \widetilde{\widetilde{C}}(x, 0) = \widetilde{\widetilde{C}}_0(x)= \max \,\left \lbrace e^{\frac{(k+1)\,x}{2}}-e^{\frac{(k-1)\,x}{2}},0\right \rbrace
\end{equation}

\noindent The classical analytical solution is then given by:

\begin{equation}\label{Classical}
    \widetilde{\widetilde{C}}_{classical} (x, \tau) =\displaystyle \frac{1}{2\,\sqrt{ \pi\,\tau}}\,\displaystyle\int_{-\infty}^{\infty}{ \widetilde{\widetilde{C}}_0 (y)\exp{\left[-\frac{(x - y)^2}{4\,\tau}\right]}}\,dy \end{equation}

\subsection{Our approach - Consistency with respect to initial conditions}

\label{Consistency}

\noindent Let us recall the natural scaling of the heat equation (\ref{Heat}). If we denote by $\widetilde{\widetilde{C}}$ a solution, then, for any strictly positive real number $\lambda$, the map

$$ (t,x) \mapsto \widetilde{\widetilde{C}}_\lambda (t,x) = \lambda \,\widetilde{\widetilde{C}}(\lambda^2t,\lambda x)$$

\noindent is also solution.\\

\noindent Following J.Y. Chemin and Cl. David \cite{CheminEtClaire1}, \cite{CheminEtClaire2}, we define the mapping, denoted by~$\cF$, from~\mbox{$L^2_{loc}(\R ) \times \R_+^\star \times \N^\star $}, by:
$$
\cF  \left ( \widetilde{\widetilde{C}}_0, \lambda, N_0 \right) = \widetilde{\widetilde{C}}_0  + \epsilon\,\sum_{j=1}^{N_0}   \lambda^{-j}\,\widetilde{\widetilde{C}}_0 \left ( \lambda^{-j}\,\cdot\right)\quad, \quad \epsilon\,\in\,\left \lbrace -1,1 \right\rbrace, \,  N_0 \,\in\,\N^\star
$$

\noindent which takes its origin in the so-called "profile decomposition theory" initiated by P. Gerard and H. Bahouri \cite{bahourigerard}. It relies on the idea that two solutions  of an evolution equations with scales that are different enough almost do not interact.\\

\noindent An important question one may ask is wether our approach does not affect the required initial conditions (\ref{CondIn2}).\\

\noindent The difficulty, here, lays in the fact that the functions at stake are not integrable on $\R$. This is the reason why we will work on $L^2_{loc}(\R)$, more precisely, on $L^2([0,S_0])$, for $S_0 \geq 0$.

\noindent The main property of~$\cF$ is given by the proposition that follows.
\begin{proposition}
\label{L2FunctionLemma}
{\sl  There exists a strictly positive constant $\lambda_0$ such that:
\begin{equation}
\label{L2FunctionLemmaeq1}
\lambda \geq \lambda_0 \quad \Rightarrow \quad
   \left \|  \cF  \left ( \widetilde{\widetilde{C}}_0, \lambda, N_0 \right)- \widetilde{\widetilde{C}}_0  \right\|^2_{  L^2(  [0,S_0] ) }= o(1)
\end{equation}
}
\end{proposition}
\vskip 1cm

\begin{proof}

One has:

$$\begin{array}{ccc}\left \|  \cF  \left ( \widetilde{\widetilde{C}}_0, \lambda, N_0 \right)- \widetilde{\widetilde{C}}_0  \right\|^2_{  L^2(  [0,S_0] ) }&=&
\left \| \displaystyle \sum_{j=1}^{N_0} \lambda^{- j}\, \widetilde{\widetilde{C}}_0 \left ( \lambda^{-j}\,\cdot\right) \right\|^2_{  L^2( [0,S_0] ) }
\\
&=&   \displaystyle \sum_{j=1}^{N_0} \lambda^{-2j}\,\left \|\widetilde{\widetilde{C}}_0 \left ( \lambda^{-j}\,\cdot\right) \right\|^2_{  L^2( [0,S_0] ) }\\
&&+2\,\displaystyle \sum_{ 1 \leq j<k \leq N_0} \lambda^{- j-k}\,\left (\widetilde{\widetilde{C}}_0 \left ( \lambda^{-j}\,\cdot\right),
\widetilde{\widetilde{C}}_0 \left ( \lambda^{-k}\,\cdot\right) \right)_{  L^2( [0,S_0] )}\\
&\leq& \left (\displaystyle \sum_{j=1}^{N_0} \lambda^{-  j} \right)\,\left \|\widetilde{\widetilde{C}}_0   \right\|^2_{  L^2( [0,S_0] ) }\\
&&+2\,\displaystyle \sum_{  1 \leq j<k \leq N_0} \lambda^{- j-k}\,\left (\widetilde{\widetilde{C}}_0 \left ( \lambda^{-j}\right) \,\cdot\widetilde{\widetilde{C}}_0 \left ( \lambda^{-k}\,\right) \right)_{  L^2(  [0,S_0]  )}\\
&=&  \lambda^{-1}\,\displaystyle \frac{1-\lambda^{-N_0}}{1-\lambda^{-1}}\,\left \|\widetilde{\widetilde{C}}_0   \right\|^2_{  L^2( [0,S_0] ) }\\
&& +2\,\displaystyle \sum_{   1 \leq j< k \leq N_0} \lambda^{- j-k}\,\left (\widetilde{\widetilde{C}}_0 \left ( \lambda^{-j}\,\cdot\right) ,\widetilde{\widetilde{C}}_0 \left ( \lambda^{-k}\,\cdot\right) \right)_{  L^2(  [0,S_0]  )}\\
\end{array}
$$

\noindent due to:

$$\begin{array}{ccc}  \lambda^{- 2j}\,\left \|\widetilde{\widetilde{C}}_0 \left ( \lambda^{-j}\,\cdot\right) \right\|^2_{  L^2( [0,S_0] ) }
&=&
\lambda^{- 2j}\,\displaystyle \int_{ 0 }^{S_0} \widetilde{\widetilde{C}}_0^2 \left ( \lambda^{-j}\,\cdot\right) \,dS \\
&=&
\lambda^{- 2j}\,\displaystyle \int_{ 0 }^{ \lambda^{-j}\,S_0} \widetilde{\widetilde{C}}_0^2 \left (  \cdot\right) \,\lambda^{ j}\,dS \\
& \leq &
\lambda^{-  j}\, \displaystyle \int_{ 0 }^{  S_0} \widetilde{\widetilde{C}}_0^2 \left (  \cdot\right)  \,dS \\
 &=& \lambda^{-  j}\,\left \|\widetilde{\widetilde{C}}_0  \right\|^2_{  L^2( [0,S_0] ) }
\end{array}
$$

\noindent and, for any set of integers \mbox{$(j,k)\,\in\,\left \lbrace 1, \hdots, N_0 \right \rbrace$} such that $j< k$, provided that the scaling factor $\lambda$ is great enough, the pseudo-orthogonality, or the fact that scales that are different enough almost do not interact, yields:

$$\begin{array}{ccc}  \lambda^{- j-k}\,\left (\widetilde{\widetilde{C}}_0 \left ( \lambda^{-j}\,\cdot\right),\widetilde{\widetilde{C}}_0 \left ( \lambda^{-k}\,\cdot\right) \right)_{  L^2(  [0,S_0]  )}
&=&
\lambda^{ -k}\,\displaystyle \int_{ 0 }^{\lambda^{-  j}\,S_0} \widetilde{\widetilde{C}}_0  \left ( \cdot\right) \,\widetilde{\widetilde{C}}_0  \left ( \lambda^{-(k-j)}\,\cdot\right) \,dS \\
&=& o(1)
\end{array}
$$

\noindent For any strictly positive real number $\varepsilon$, one easily can find the threshold value $\lambda_0$ such that:

$$ \forall \,\lambda> \lambda_0 \, : \quad\lambda^{-1}\,\displaystyle \frac{1-\lambda^{-N_0}}{1-\lambda^{-1}}\,\left \|\widetilde{\widetilde{C}}_0   \right\|^2_{  L^2( [0,S_0] ) }  \leq
 \varepsilon   $$
\end{proof}
\vskip 1cm

\begin{remark}

The mapping~$\cal F$ shows that the control depends on~$N_0$ and~$\lambda$. As~$\lambda$ increases, one will require smaller values of the integer~$N_0$:

$$ N_0 \ln \lambda   \leq \text{exp} \left \lbrace -\displaystyle \frac{1}{N_0}\, \ln \left (  1-
 \varepsilon \, \lambda\, \displaystyle \frac{  1-\lambda^{-1}}{ \,\,\,\,\,\left \|\widetilde{\widetilde{C}}_0   \right\|^2_{  L^2( [0,S_0] ) } }  \right ) \right \rbrace
   $$

\end{remark}

\vskip 1cm

\noindent The illustration of the above theoretical results can be seen through the following numerical results. We hereafter display the variations of:

\begin{enumerate}
\item[$\rightsquigarrow$] the initial condition function $ \widetilde{\widetilde{C}}_0$ given by (\ref{CondIn2}), in red ;

\item[$\rightsquigarrow$] the initial condition function $\widetilde{\widetilde{C}}_{ 0,\lambda}$ given by:

$$ \widetilde{\widetilde{C}}_{ 0,\lambda}(x)=\widetilde{\widetilde{C}}_{0}(x)+\epsilon\,
\displaystyle \sum_{j=1}^{N_0} \widetilde{\widetilde{C}}_{0, \lambda,j }(x )=\widetilde{\widetilde{C}}_{0}(x)+\epsilon\,\displaystyle \sum_{j=1}^{N_0}\displaystyle \frac{1}{\lambda^j}\,\widetilde{\widetilde{C}}_0\left ( \displaystyle \frac{x}{\lambda^j} \right )
=\widetilde{\widetilde{C}}_{0}(x)+\epsilon\,\displaystyle \sum_{j=1}^{N_0}\displaystyle \frac{1}{\lambda^j}\,\max \,\left \lbrace e^{\frac{(k+1)\,x}{2\,\lambda^j}}-e^{\frac{(k-1)\,x}{2\,\lambda^j}},0\right \rbrace  $$

 in green ;
\end{enumerate}

\noindent in the case where:

$$ r = 0.06 \quad , \quad \sigma = 0.3 \quad , \quad  E = 100$$

\noindent If one chooses values of $\lambda$ greater or equal to 10, the required initial conditions (\ref{CondIn2}) are not affected.

\begin{figure}[h!]
 \center{\psfig{height=5cm,width=7cm,angle=0,file=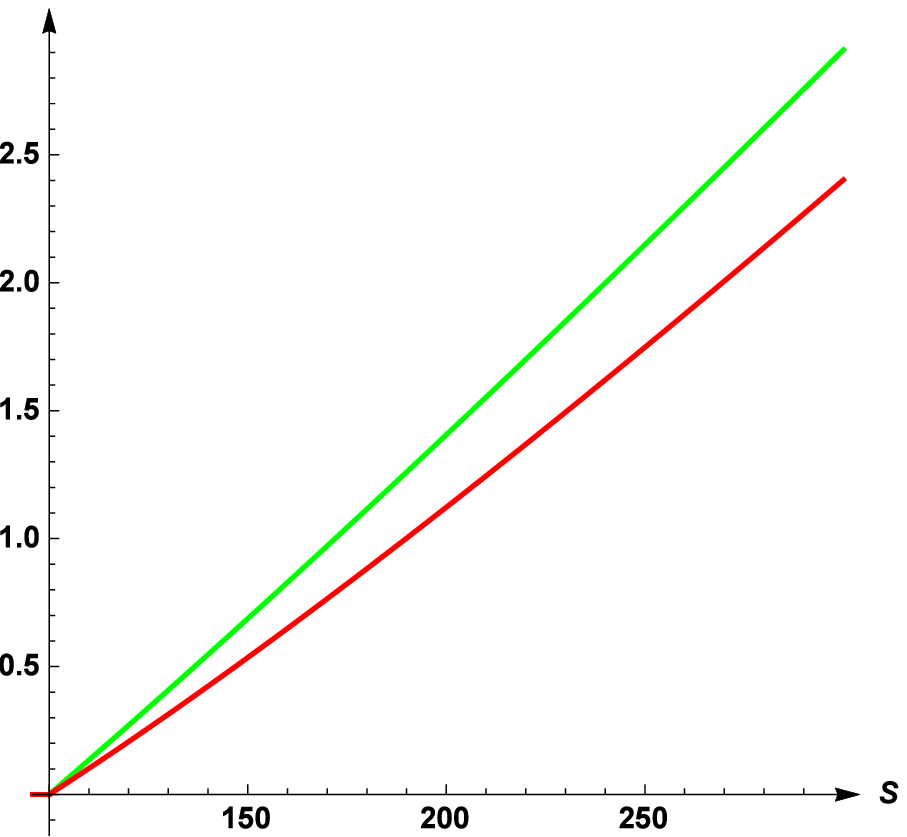}} \\
 \caption{The graph of $ \widetilde{\widetilde{C}}_0$ and $\widetilde{\widetilde{C}}_{ 0,\lambda}$ as functions of stock price $S$, for $\lambda=2$.}
  \end{figure}

\begin{figure}[h!]
 \center{\psfig{height=5cm,width=7cm,angle=0,file=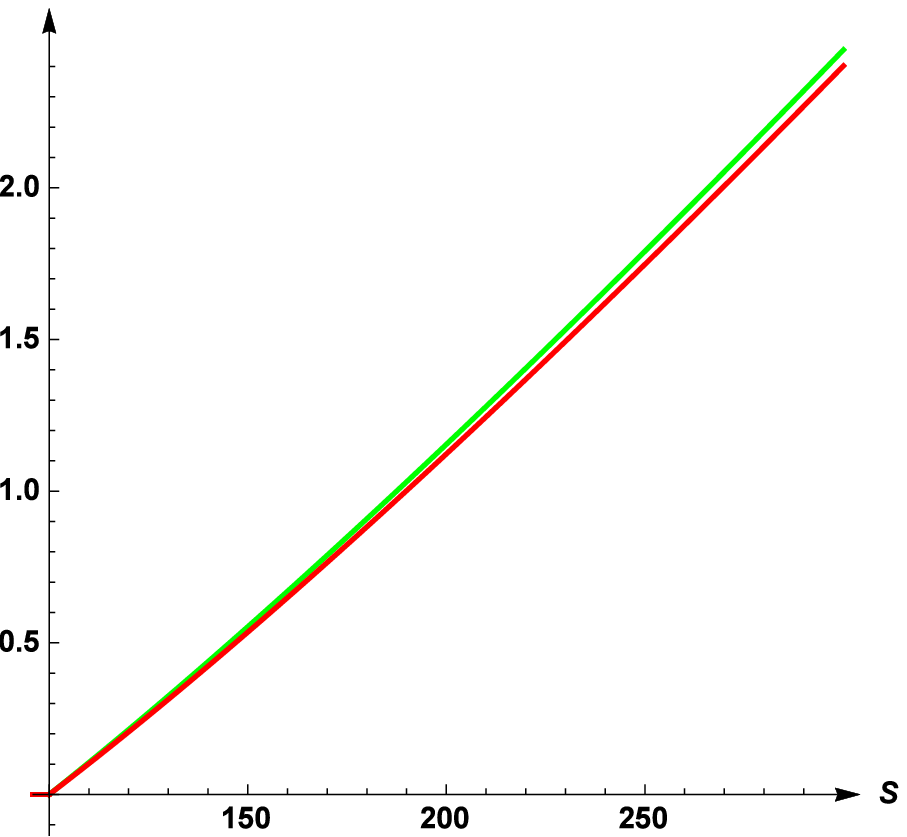}} \\
 \caption{The graph of $ \widetilde{\widetilde{C}}_0$ and $\widetilde{\widetilde{C}}_{ 0,\lambda}$ as functions of stock price $S$, for $\lambda=5$.}
  \end{figure}

\begin{figure}[h!]
 \center{\psfig{height=5cm,width=7cm,angle=0,file=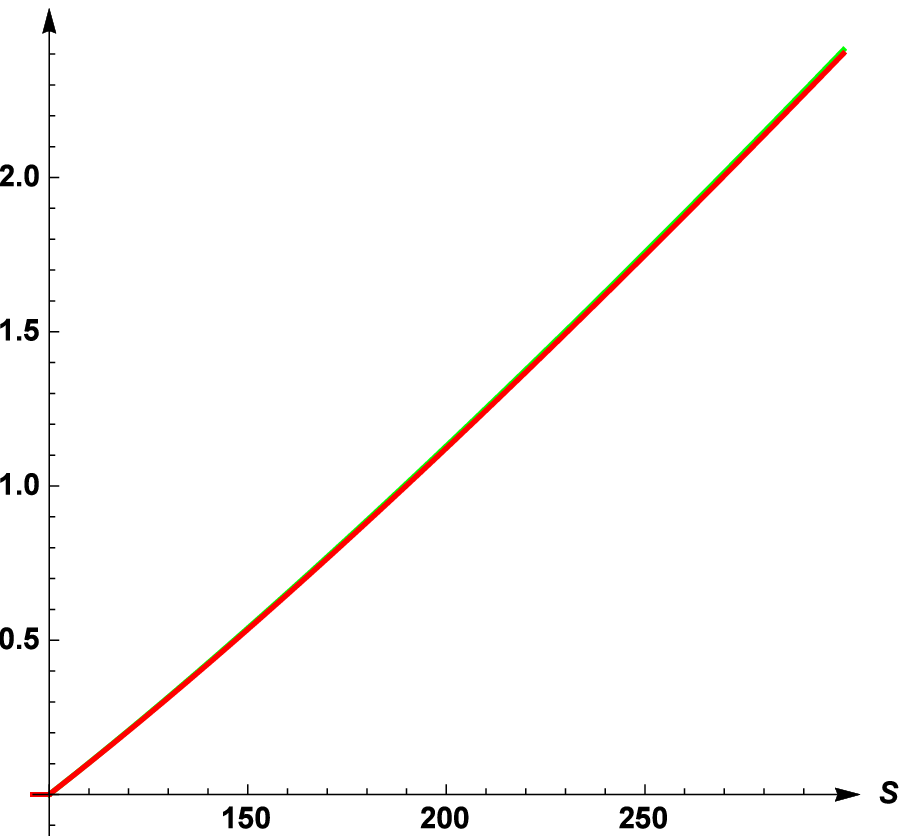}} \\
 \caption{The graph of $ \widetilde{\widetilde{C}}_0$ and $\widetilde{\widetilde{C}}_{ 0,\lambda}$ as functions of stock price $S$, for $\lambda=10$.}
  \end{figure}

\begin{figure}[h!]
 \center{\psfig{height=5cm,width=7cm,angle=0,file=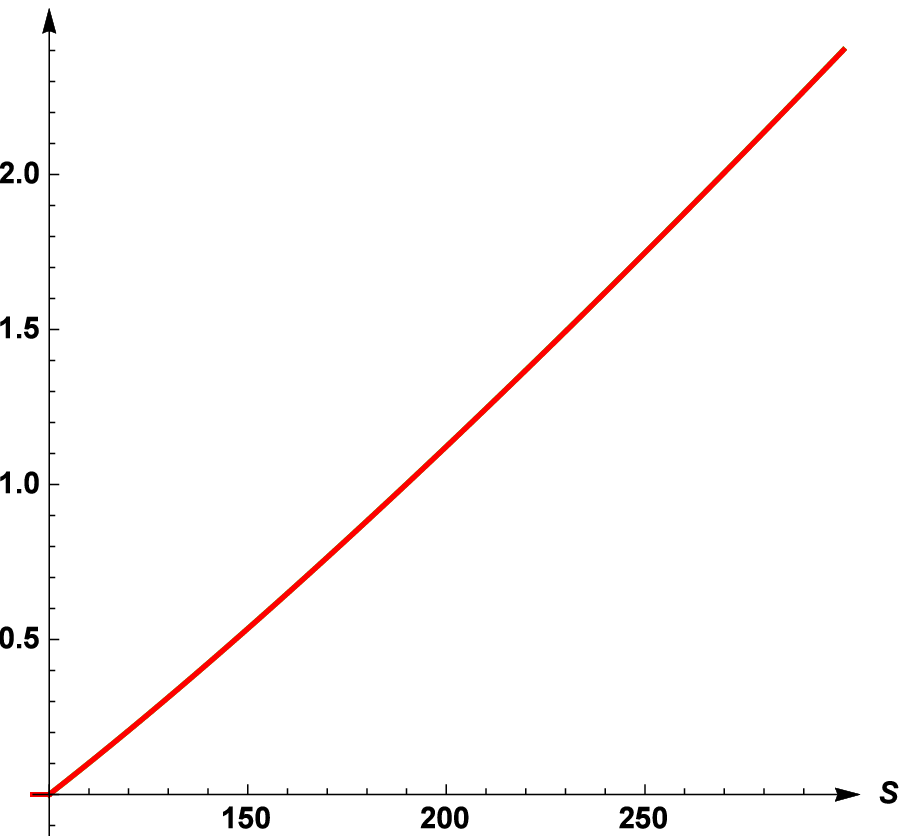}} \\
 \caption{The graph of $ \widetilde{\widetilde{C}}_0$ and $\widetilde{\widetilde{C}}_{ 0,\lambda}$ as functions of stock price $S$, for $\lambda=50$.}
  \end{figure}

\begin{figure}[h!]
 \center{\psfig{height=5cm,width=7cm,angle=0,file=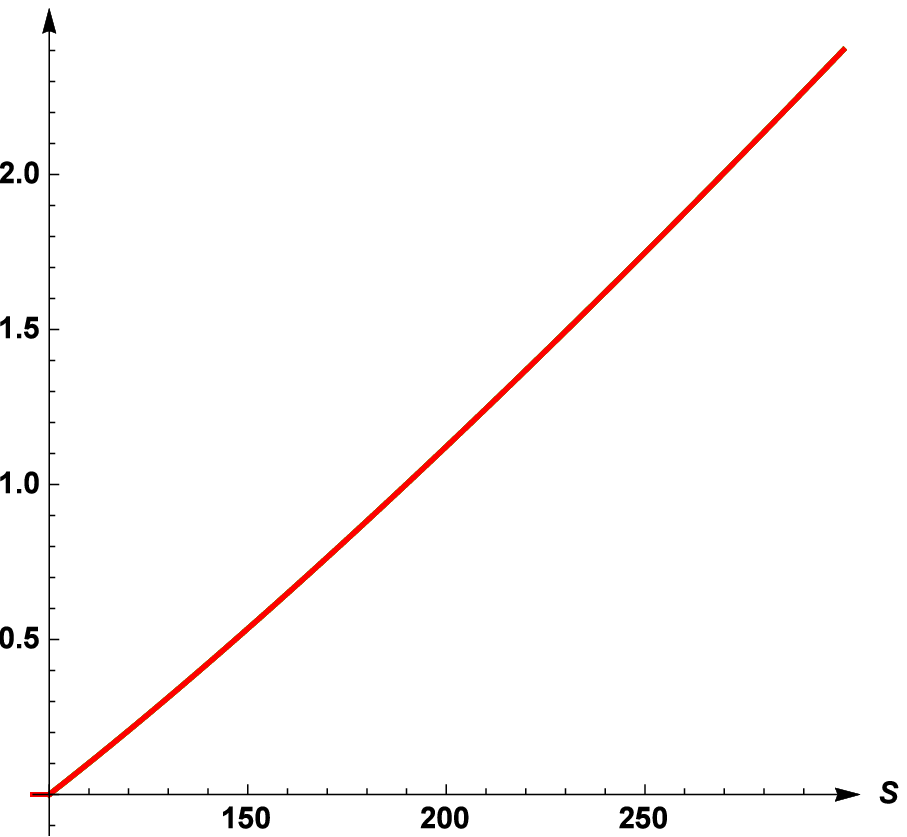}} \\
 \caption{The graph of $ \widetilde{\widetilde{C}}_0$ and $\widetilde{\widetilde{C}}_{ 0,\lambda}$ as functions of stock price $S$, for $\lambda=100$.}
  \end{figure}

\subsection{Analytic results}
\noindent As previously, we consider, for $\Lambda>\Lambda_0$, initial data of the form:

$$ \widetilde{\widetilde{C}}_{ 0,\lambda}(x)=\widetilde{\widetilde{C}}_{0}(x)+\epsilon\,\displaystyle \sum_{j=1}^{N_0} \widetilde{\widetilde{C}}_{0, \lambda,j }(x )=\widetilde{\widetilde{C}}_{0}(x)+\epsilon\,\displaystyle \sum_{j=1}^{N_0}\displaystyle \frac{1}{\lambda^j}\,\widetilde{\widetilde{C}}_0\left ( \displaystyle \frac{x}{\lambda^j} \right )
=\widetilde{\widetilde{C}}_{0}(x)+\epsilon\,\displaystyle \sum_{j=1}^{N_0}\displaystyle \frac{1}{\lambda^j}\,\max \,\left \lbrace e^{\frac{(k+1)\,x}{2\,\lambda^j}}-e^{\frac{(k-1)\,x}{2\,\lambda^j}},0\right \rbrace  $$

\noindent The related exact analytic solution $\widetilde{\widetilde{C}}$, which is a function of $x$, $\tau$, and of the scaling parameter $\lambda$, is given by:

$$\widetilde{\widetilde{C}}(x,\tau,\lambda)=\widetilde{\widetilde{C}}_{classical} (x, \tau)+\epsilon\,\displaystyle \frac{1}{2\,\sqrt{\pi\,\tau}}\,\displaystyle \int_{-\infty}^{+\infty}
 \displaystyle \sum_{j=1}^{N_0} \widetilde{\widetilde{C}}_{0, \lambda,j }(y) \,e^{-\frac{(x-y)^2}{4\,\tau}}\,dy$$

\noindent For any integer $j$ in \mbox{$ \left \lbrace 1, \hdots, N_0 \right \rbrace$ }:

$$\begin{array}{ccc}\displaystyle \frac{1}{2\,\sqrt{\pi\,\tau}}\, \displaystyle \int_{-\infty}^{+\infty}
  \widetilde{\widetilde{C}}_{0, \lambda,j }(y) \,e^{-\frac{(x-y)^2}{4\,\tau}}\,dy &=&\displaystyle \frac{1}{ \sqrt{\pi }\,\lambda^j}\, \displaystyle \int_{- \infty}^{+\infty}  \max \,\left \lbrace e^{\frac{(k+1)\,x}{2\,\lambda^j}}-e^{\frac{(k-1)\,x}{2\,\lambda^j}},0\right \rbrace   \,e^{- z^2 }\,dz \\ \\
   &=&\displaystyle \frac{1}{ \sqrt{\pi }\,\lambda^j}\, \displaystyle \int_{-\frac{x}{ 2\,\sqrt{ \tau} }}^{+\infty}  e^{ \frac{ (k+1)\,(x+2\,z\,\sqrt{ \tau})}{2\,\lambda^j} } \,e^{- z^2 }\,dz \\
   &&- \displaystyle \frac{1}{ \sqrt{\pi }\,\lambda^j}\, \displaystyle \int_{ -\frac{x}{ 2\,\sqrt{ \tau} }}^{+\infty}   e^{\frac{  (k-1)\,(x+2\,z\,\sqrt{ \tau}) }{2\,\lambda^j}} \,e^{- z^2 }\,dz \\ \\
&=&\ \displaystyle \frac{1}{ \sqrt{\pi }\,\lambda^j}\, \displaystyle \int_{-\frac{x}{ 2\,\sqrt{ \tau} }}^{+\infty}  e^{ \frac{ (k+1)\,x  }{2\,\lambda^j }} \,e^{- \left (z+\frac{(k+1)\,\sqrt{\tau}}{4\,\lambda^j}\right)^2 }\,dz
\\
&&- \displaystyle \frac{1}{ \sqrt{\pi }\,\lambda^j}\, \displaystyle \int_{-\frac{x}{ 2\,\sqrt{ \tau} }}^{+\infty}  e^{ (k-1)\,x +(k-1)^2\,\tau} \,e^{- \left (z+\frac{(k-1)\,\sqrt{\tau}}{2\, \lambda^j}\right)^2 }\,dz\\ \\
&=& \displaystyle \frac{e^{\frac{ (k+1)\,x}{2\,\lambda^j} +\frac{(k+1)^2\,\tau}{4\,\lambda^{2j}} } }{ \sqrt{\pi }\,\lambda^j}\, \displaystyle \int_{-\frac{x}{ 2\,\sqrt{ \tau} }}^{+\infty}   e^{- \left (z+\frac{(k+1)\,\sqrt{\tau}}{2\,\lambda^j}\right)^2 }\,dz
\\ &&- \displaystyle \frac{e^{ \frac{(k-1)\,x}{2\,\lambda^j} +\frac{(k-1)^2\,\tau}{4\,\lambda^{2j}}}  }{ \sqrt{\pi }\,\lambda^j}\, \displaystyle \int_{-\frac{x}{ 2\,\sqrt{ \tau} }}^{+\infty}  e^{- \left (z+\frac{(k-1)\,\sqrt{\tau}}{2\,\lambda^j}\right)^2 }\,dz\\ \\
&=& \displaystyle \frac{e^{ \frac{(k+1)\,x}{2\,\lambda^{ j}} +\frac{(k+1)^2\,\tau}{4\,\lambda^{2j}}}  }{ \sqrt{\pi }\,\lambda^j}\, \displaystyle \int_{-\frac{x}{ 2\,\sqrt{ \tau} }+ \frac{(k+1)\,\sqrt{\tau}}{2\,\lambda^j}}^{+\infty}   e^{-  z   }\,dz
\\ &&-  \displaystyle \frac{e^{ \frac{(k-1)\,x}{2\,\lambda^j} +\frac{(k-1)^2\,\tau}{4\,\lambda^{2j}}} }{ \sqrt{\pi }\,\lambda^j}\, \displaystyle \int_{-\frac{x}{ 2\,\sqrt{ \tau} }+ \frac{(k -1)\,\sqrt{\tau}}{2\,\lambda^j}}^{+\infty}  e^{- z^2 }\,dz\\ \\
&=& \displaystyle \frac{e^{ \frac{(k+1)\,x}{2\,\lambda^j} +\frac{(k+1)^2\,\tau}{4\,\lambda^{2j}}}   }{ \sqrt{\pi }\,\lambda^j}\,\displaystyle \frac {\sqrt{\pi}}{2}\,Erf_c\left (-\frac{x}{ 2\,\sqrt{ \tau} }+  \frac{(k+1)\,\sqrt{\tau}}{2\,\lambda^j} \right) \\
  &&- \displaystyle \frac{e^{ \frac{(k-1)\,x}{2\,\lambda^{ j}} +\frac{(k-1)^2\,\tau}{4\,\lambda^{2j}}} }{ \sqrt{\pi }\,\lambda^j}\,\displaystyle \frac {\sqrt{\pi}}{2}\,Erf_c\left (-\frac{x}{ 2\,\sqrt{ \tau} }+ \frac{ (k-1)\,\sqrt{\tau}}{2\,\lambda^j} \right) \\ \\
&=& \displaystyle \frac{1}{\lambda^j}\,\displaystyle  e^{ \frac{(k+1)\,x}{2\,\lambda^j} +\frac{(k+1)^2\,\tau}{4\,\lambda^{2j}}}    \,N\left ( \frac{\sqrt{2}\,x}{ 2\,\sqrt{ \tau} }- \sqrt{2}\, \frac{(k+1)\,\sqrt{\tau}}{2\,\lambda^j} \right)  \\
 &&-  \displaystyle  \frac{1}{ \lambda^j} \,e^{ \frac{(k-1)\,x}{2\,\lambda^j} +\frac{(k-1)^2\,\tau}{4\,\lambda^{2j}}}  \, N\left (-\frac{\sqrt{2}\,x}{ 2\,\sqrt{ \tau} }- \sqrt{2}\, \frac{(k-1)\,\sqrt{\tau}}{2\,\lambda^j} \right) \\ \\
\end{array}$$

\noindent where the complementary error function $Erf_c$ is defined, for any real number $x$, by:

$$
  Erf_c(x)   = \displaystyle \frac{2}{\sqrt{\pi}} \, \displaystyle\int_x^{+\infty} e^{-t^2}\,  dt  $$

\noindent while the normal (gaussian) cumulative distribution function is given, for any real number $d$, by:

$$
  N(d)   = \displaystyle \frac{1}{\sqrt{2\pi}} \, \displaystyle\int_{ -\infty}^{d} e^{-\frac{t^2}{2}}\,  dt
    =\displaystyle \frac{1}{\sqrt{2\pi}}\, \sqrt{2}\, \displaystyle\int_{-\frac{d}{\sqrt{2}}}^{ +\infty}e^{-  t^2 }\,  dt
    =\displaystyle \frac{1}{ 2}\, Erf_c\left (-\displaystyle \frac{d}{\sqrt{2}}\right)
     $$

\noindent Thus:

$$\begin{array}{ccc}
\widetilde{\widetilde{C}}(x,\tau,\lambda  )&= &\widetilde{\widetilde{C}}_{classical} (x, \tau)+\epsilon\,\displaystyle \sum_{j=1}^{N_0}\displaystyle \frac{e^{ \frac{(k+1)\,x}{2\,\lambda^j} +\frac{(k+1)^2\,\tau}{4\,\lambda^{2j}}}   }{ 2\,\lambda^j} \,Erf_c\left (-\frac{x}{ 2\,\sqrt{ \tau} }+  \frac{(k+1)\,\sqrt{\tau}}{2\,\lambda^j} \right) \\
  &&-\epsilon\,\displaystyle \sum_{j=1}^{N_1}\displaystyle \frac{e^{ \frac{(k-1)\,x}{2\,\lambda^{ j}} +\frac{(k-1)^2\,\tau}{4\,\lambda^{2j}}} }{ 2\,\lambda^j}\, \,Erf_c\left (-\frac{x}{ 2\,\sqrt{ \tau} }+ \frac{ (k-1)\,\sqrt{\tau}}{2\,\lambda^j} \right) \\ \\
\end{array}$$

\noindent It is interesting to note that:

$$\begin{array}{ccc}
\widetilde{\widetilde{C}}(x,\tau,\lambda  )&= &\widetilde{\widetilde{C}}_{classical} (x, \tau)+
\epsilon\,\displaystyle \sum_{j=1}^{N_0} \displaystyle \frac{1}{\lambda^j}\,\widetilde{\widetilde{C}}_{classical} \left (\displaystyle \frac{x}{\lambda^j}, \displaystyle \frac{\tau}{\lambda^{2j}}\right)
\end{array}$$

\noindent and to notice, thus, that it includes, as expected, the natural scaling of the heat equation (\ref{Heat}).\\

\noindent One easily goes back to the call function $C$ through

$$C(S,t)=E\,\widetilde{C} (x,\tau)= E\,e^{-\alpha\,x-\beta\,\tau}\,\widetilde{\widetilde{C}}(x,\tau,\lambda )\quad \text{with} \quad x =\ln \displaystyle \frac{S}{E}$$

\section{Results}

\label{Results}

\noindent In finance, the sensitivity of a portfolio to changes in parameters values can be measured through what commonly call "the Greeks", i.e.:

\begin{enumerate}

\item[\emph{i}.] the Delta \mbox{$\Delta=\displaystyle \frac{\partial C }{\partial S} \,\in\,[0,1]$}, which enables one to quantify the risk, and is thus the most important Greek.

\item[\emph{ii}.] The Gamma $\Gamma=\displaystyle \frac{\partial^2 C }{\partial  S^2} \geq 0$.

\item[\emph{iii}.] The Vega (the name of which comes from the form of the greek letter $\nu$) $\nu=\displaystyle \frac{\partial  C }{\partial \sigma}$.

\item[\emph{iv}.] The Theta $\Theta=\displaystyle \frac{\partial  C }{\partial t}$.

\item[\emph{v}.] The rho $\rho=\displaystyle \frac{\partial  C }{\partial r}$.\\

\end{enumerate}

\noindent   The good strategy, for traders, is to have delta-neutral positions at least once a day, and, whenever the opportunity arises, to improve the Gamma and the Vega
\cite {BlackScholes}.

\subsection{Control of the Delta and Gamma}

\noindent To test our approach, we have choosen to compare, first:

\begin{enumerate}
\item[$\rightsquigarrow$] the classical Delta $\Delta_{classical}$ and Gamma $\Gamma_{classical}$ ;

\item[$\rightsquigarrow$] the ones of our approach.
\end{enumerate}

\noindent Due to the decomposition
$$\begin{array}{ccc}
\widetilde{\widetilde{C}}(x,\tau,\lambda  )&= &\widetilde{\widetilde{C}}_{classical} (x, \tau)+
\epsilon\,\displaystyle \sum_{j=1}^{N_0} \displaystyle \frac{1}{\lambda^j}\,\widetilde{\widetilde{C}}_{classical} \left (\displaystyle \frac{x}{\lambda^j}, \displaystyle \frac{\tau}{\lambda^{2j}}\right)
\end{array}$$

\noindent the change of variables $x =\ln \displaystyle \frac{S}{E}$ leads to:

\begin{equation}\label{Delta}\begin{array}{ccc}\displaystyle \frac{\partial C }{\partial S}
 &= & \widetilde{\widetilde{\Delta}}_{classical} (x, \tau)+
\epsilon\, \displaystyle \frac{E}{S} \,\displaystyle \sum_{j=1}^{N_0} \displaystyle \frac{1}{\lambda^j}\,\displaystyle\frac{\partial   }{\partial x} \left [\widetilde{\widetilde{C}}_{classical} \left (\displaystyle \frac{x}{\lambda^j}, \displaystyle \frac{\tau}{\lambda^{2j}}\right) \right]\\
&= &\widetilde{\widetilde{\Delta}}_{classical}   (x, \tau)+
\epsilon\, \displaystyle \frac{E}{S} \, \displaystyle \sum_{j=1}^{N_0} \, \displaystyle \frac{1}{\lambda^j}\,\displaystyle\frac{\partial   }{\partial x} \left [
   \widetilde{\widetilde{C}}_{classical} \left (\displaystyle \frac{x}{\lambda^j}, \displaystyle \frac{\tau}{\lambda^{2j}}\right) \right]\\
   &= &\widetilde{\widetilde{\Delta}}_{classical} (x, \tau)+
\epsilon\,   \displaystyle \sum_{j=1}^{N_0} \, \displaystyle \frac{1}{\lambda^{2\,j}}\,
   \widetilde{\widetilde{\Delta}}_{classical}  \left (\displaystyle \frac{x}{\lambda^j}, \displaystyle \frac{\tau}{\lambda^{2j}}\right) \\
 \end{array}
 \end{equation}

\noindent where we have set:

$$\widetilde{\widetilde{\Delta}}_{classical}  (x, \tau)=\Delta_{classical}(S,t)=\displaystyle \frac{\partial C }{\partial S} $$

\noindent It appears thus that:

\begin{enumerate}
\item[$\rightsquigarrow$] for $\epsilon=1$, one can increase the Delta;

\item[$\rightsquigarrow$] for $\epsilon=-1$, one can decrease the Delta.
\end{enumerate}

\noindent In practice, it seems interesting to determine a suitable value $\lambda_0$ of the shape parameter $\lambda$ such that:

$$\Delta=\left(1+\epsilon\,\eta_0 \right)\,\Delta_{classical}  \quad ,\quad \eta_0\,\in\,\left ]0,1\right [$$

\noindent It can be achieved through a series expansion of the quantity \mbox{$\displaystyle \frac{1}{\lambda^{2\,j}}\,
   \widetilde{\widetilde{\Delta}}_{classical} \left (\displaystyle \frac{x}{\lambda^j}, \displaystyle \frac{\tau}{\lambda^{2j}}\right)$}.\\

\noindent One has:

$$\begin{array}{ccc}
\widetilde{\widetilde{\Delta}}_{classical} (x,\tau) &= &
 \displaystyle \frac{E}{S} \, \displaystyle\frac{\partial   }{\partial x} \left [\widetilde{\widetilde{C}}(x,\tau   ) \right]\\
 &= &
 \displaystyle \frac{E}{S} \, \displaystyle\frac{\partial   }{\partial x} \left [ \displaystyle \frac{e^{ \frac{(k+1)\,x}{2} +\frac{(k+1)^2\,\tau}{4}}   }{ \sqrt{\pi }}\,\displaystyle \frac {\sqrt{\pi}}{2}\,Erf_c\left (-\frac{x}{ 2\,\sqrt{ \tau} }+  \frac{(k+1)\,\sqrt{\tau}}{2} \right)     \right]\\
&  & -
 \displaystyle \frac{E}{S} \, \displaystyle\frac{\partial   }{\partial x} \left [ \displaystyle \frac{e^{ \frac{(k-1)\,x}{2} +\frac{(k-1)^2\,\tau}{4}} }{ \sqrt{\pi }}\,\displaystyle \frac {\sqrt{\pi}}{2}\,Erf_c\left (-\frac{x}{ 2\,\sqrt{ \tau} }+ \frac{ (k-1)\,\sqrt{\tau}}{2} \right) \right]\\
&= &
 \displaystyle \frac{E}{S} \, \displaystyle\frac{(k+1) }{2} \, \displaystyle \frac{e^{ \frac{(k+1)\,x}{2} +\frac{(k+1)^2\,\tau}{4}}   }{ \sqrt{\pi }}\,\displaystyle \frac {\sqrt{\pi}}{2}\,Erf_c\left (-\frac{x}{ 2\,\sqrt{ \tau} }+  \frac{(k+1)\,\sqrt{\tau}}{2} \right)     \\
&&+ \displaystyle\frac{1}{ 2\,\sqrt{ \tau} }\,
 \displaystyle \frac{E}{S} \,   \displaystyle \frac{e^{ \frac{(k+1)\,x}{2} +\frac{(k+1)^2\,\tau}{4}}   }{ \sqrt{\pi }}\,
    e^{ \left (-\frac{x}{ 2\,\sqrt{ \tau} }+ \frac{ (k-1)\,\sqrt{\tau}}{2} \right)^2}   \\
&  & -
 \displaystyle \frac{E}{S} \, \displaystyle \frac{(k-1) }{2} \, \displaystyle \frac{e^{ \frac{(k-1)\,x}{2} +\frac{(k-1)^2\,\tau}{4}} }{ \sqrt{\pi }}\,\displaystyle \frac {\sqrt{\pi}}{2}\,Erf_c\left (-\frac{x}{ 2\,\sqrt{ \tau} }+ \frac{ (k-1)\,\sqrt{\tau}}{2} \right)  \\
&  & -\displaystyle\frac{1}{ 2\,\sqrt{ \tau} }\,
 \displaystyle \frac{E}{S} \,   \displaystyle \frac{e^{ \frac{(k-1)\,x}{2} +\frac{(k-1)^2\,\tau}{4}} }{ \sqrt{\pi }}\,
   \,e^{ \left (-\frac{x}{ 2\,\sqrt{ \tau} }+ \frac{ (k-1)\,\sqrt{\tau}}{2} \right)^2}  \\
\end{array}
$$

\noindent and, therefore:

\begin{equation}
\label{DeltaClassiqueShape}\begin{array}{ccc}
\widetilde{\widetilde{\Delta}} \left (\displaystyle \frac{x}{\lambda^j}, \displaystyle \frac{\tau}{\lambda^{2j}}\right)
&= &
 \displaystyle \frac{E}{S} \, \displaystyle\frac{(k+1) }{2} \, \displaystyle \frac{e^{ \frac{(k+1)\,x}{2\lambda^j\,} +\frac{(k+1)^2\,\tau}{4\,\lambda^j}}   }{ \sqrt{\pi }}\,\displaystyle \frac {\sqrt{\pi}}{2}\,Erf_c\left (-\frac{x}{ 2\,\lambda^j\,\sqrt{ \tau} }+  \frac{(k+1)\,\sqrt{\tau}}{2\,\lambda^j} \right)     \\
&&+ \displaystyle\frac{\lambda^j}{ 2\,\sqrt{ \tau} }\,
 \displaystyle \frac{E}{S} \,   \displaystyle \frac{e^{ \frac{(k+1)\,x}{2\lambda^j\,} +\frac{(k+1)^2\,\tau}{4\,\lambda^j }}   }{ \sqrt{\pi }}\,  e^{ \left (-\frac{x}{ 2\,\lambda^j\,\sqrt{ \tau} }+ \frac{ (k-1)\,\sqrt{\tau}}{2\,\lambda^j } \right)^2}   \\
&  & -
 \displaystyle \frac{E}{S} \, \displaystyle \frac{(k-1) }{2} \, \displaystyle \frac{e^{ \frac{(k-1)\,x}{2\,\lambda^j } +\frac{(k-1)^2\,\tau}{4\,\lambda^j}} }{ \sqrt{\pi }}\,\displaystyle \frac {\sqrt{\pi}}{2}\,Erf_c\left (-\frac{x}{ 2\,\lambda^j\,\sqrt{ \tau} }+ \frac{ (k-1)\,\sqrt{\tau}}{2\,\lambda^j} \right)  \\
&  & -\displaystyle\frac{\lambda^j}{ 2\,\sqrt{ \tau} }\,
 \displaystyle \frac{E}{S} \,   \displaystyle \frac{e^{ \frac{(k-1)\,x}{2\lambda^j\,} +\frac{(k-1)^2\,\tau}{4\,\lambda^j }} }{ \sqrt{\pi }}\, e^{ \left (-\frac{x}{ 2\,\lambda^j\,\sqrt{ \tau} }+ \frac{ (k-1)\,\sqrt{\tau}}{2\,\lambda^j } \right)^2}  \\
\end{array}
\end{equation}

\noindent One requires then the series expansion of the $Erf_c$ function in 0, which is given, for $z \,\in\,\R$, by:

$$\begin{array}{ccc}
Erf_c (z)&=& \displaystyle \frac{2}{\sqrt{\pi}}\, \displaystyle \sum_{n=0}^{+\infty} \displaystyle \frac{(-1)^n \,z^{2n+1}}{ (2n+1)\,n\,!}
\end{array}$$

\noindent For any integer $j$ in \mbox{$\left \lbrace 1,\hdots, N_0 \right \rbrace$}, a series expansion of the term \mbox{$\displaystyle \frac {2\,S \, \widetilde{\widetilde{\Delta}} \left (\displaystyle \frac{x}{\lambda^j}, \displaystyle \frac{\tau}{\lambda^{2j}}\right)}{\epsilon\,E\,\sqrt{\pi}}   $} is:

$$\begin{array}{ccc}
&&      \displaystyle \frac{1}{\lambda^{2\,j}}\, \displaystyle\frac{(k+1) }{2\,\sqrt{\pi }} \,
   \left ( \displaystyle \sum_{n=0}^{+\infty} \frac{(k+1)^n\,x^n}{2^n\,\lambda^{j\,n}\,n\,! }   \right) \,
   \left (\displaystyle \frac{2}{\sqrt{\pi}}\, \displaystyle \sum_{n=0}^{+\infty} \displaystyle \frac{(-1)^n }{ (2n+1)\,n\,!}\,\left (-\frac{x}{ 2\,\lambda^j\,\sqrt{ \tau} }+  \frac{(k+1)\,\sqrt{\tau}}{2\,\lambda^j} \right)^{2n+1} \right)     \\ \\
&&+ \displaystyle\frac{\lambda^j}{  \sqrt{ \tau} }\,
 \left ( \displaystyle \sum_{n=0}^{+\infty} \frac{(k+1)^n\,x^n}{2^n\,\lambda^{j\,n}\, n\,!}   \right) \,
   \left (1+ \left (-\frac{x}{ 2\,\lambda^j\,\sqrt{ \tau} }+ \frac{ (k+1)\,\sqrt{\tau}}{2\,\lambda^j } \right)^2+ {\cal O} \left ( \displaystyle \frac{1}{\lambda^{4j}}\right)  \right)  \\ \\
&  & -  \displaystyle \frac{(k-1) }{2\,\sqrt{\pi}} \,\left ( \displaystyle \sum_{n=0}^{+\infty} \frac{(k-1)^n\,x^n}{2^n\,\lambda^{j\,n}\,n\,! }   \right)\,
    \left (\displaystyle \frac{2}{\sqrt{\pi}}\, \displaystyle \sum_{n=0}^{+\infty} \displaystyle \frac{(-1)^n }{ (2n+1)\,n\,!}\,\left (-\frac{x}{ 2\,\lambda^j\,\sqrt{ \tau} }+  \frac{(k-1)\,\sqrt{\tau}}{2\,\lambda^j} \right)^{2n+1} \right)  \\ \\
&  & -\displaystyle\frac{\lambda^j}{ \sqrt{ \tau} }\,
     \left ( \displaystyle \sum_{n=0}^{+\infty} \frac{(k-1)^n\,x^n}{2^n\,\lambda^{j\,n}\,n\,! }   \right) \,
     \left ( \displaystyle \sum_{n=0}^{+\infty} \frac{1 }{n \,!}\,\left (-\frac{x}{ 2\,\lambda^j\,\sqrt{ \tau} }+ \frac{ (k-1)\,\sqrt{\tau}}{2\,\lambda^j } \right)^n   \right)
   \\ \\
\end{array}$$

\noindent One obtains thus, for a given set $(x,\tau)$, or, in an equivalent way, $(S,t)$, and a given value $\eta_0$, an equation in $\lambda$ that can be solved numerically.\\

\noindent In the same way as above, the change of variables $x =\ln \displaystyle \frac{S}{E}$ leads to:

$$\begin{array}{ccc}\displaystyle \frac{\partial^2 C }{\partial S^2}
 &= &\displaystyle \frac{\partial^2 C }{\partial S^2}+
\epsilon\,\displaystyle \frac{E^2}{S^2} \, \displaystyle \sum_{j=1}^{N_0} \, \displaystyle \frac{1}{\lambda^j}\,\displaystyle\frac{\partial   }{\partial S} \left [
     \displaystyle \frac{\partial   }{\partial S }\widetilde{\widetilde{C}}_{classical} \left (\displaystyle \frac{x}{\lambda^j}, \displaystyle \frac{\tau}{\lambda^{2j}}\right) \right]\\
   &= &\widetilde{\widetilde{\Gamma}}_{classical} (x, \tau)+
\epsilon\,  \displaystyle \sum_{j=1}^{N_0} \, \displaystyle \frac{1}{\lambda^{3\,j}}\,
   \Gamma_{classical} \left (\displaystyle \frac{x}{\lambda^j}, \displaystyle \frac{\tau}{\lambda^{2j}}\right) \\
 \end{array}$$

\noindent It appears thus that:

\begin{enumerate}
\item[$\rightsquigarrow$] for $\epsilon=1$, one can increase the Gamma ;

\item[$\rightsquigarrow$] for $\epsilon=-1$, one can decrease the Gamma.
\end{enumerate}

\noindent The interesting point is that, due do the relation (\ref{Delta}):

$$\begin{array}{ccc}\displaystyle \frac{\partial   }{\partial S} \left [\Delta-\Delta_{classical} \right]  &= &
\epsilon\,   \displaystyle \sum_{j=1}^{N_0} \, \displaystyle \frac{1}{\lambda^{3\,j}}\,
   \widetilde{\widetilde{\Gamma}}_{classical}  \left (\displaystyle \frac{x}{\lambda^j}, \displaystyle \frac{\tau}{\lambda^{2j}}\right) \\
 \end{array}$$

 \noindent this expression being of the same sign as $\epsilon$. Hence:\\

\begin{enumerate}
\item[$\rightsquigarrow$] for $\epsilon=1$, the positive quantity $\Delta-\Delta_{classical} $  is an increasing function of the asset stock price $S$;

\item[$\rightsquigarrow$] for $\epsilon=-1$, the negative quantity $\Delta-\Delta_{classical} $  is a decreasing function of the asset stock price $S$.
\end{enumerate}

\noindent It appears thus that if one finds a value $\lambda_0$ of the scale parameter such that

\begin{enumerate}
\item[$\rightsquigarrow$] $\Delta=\left(1+  \eta_0 \right)\,\Delta_{classical}  \quad ,\quad \eta_0\,\in\,\left ]0,1\right [$, one will increase the Delta on the study interval ;

\item[$\rightsquigarrow$] $\Delta=\left(1-  \eta_0 \right)\,\Delta_{classical}  \quad ,\quad \eta_0\,\in\,\left ]0,1\right [$, one will decrease the Delta on the study interval.

\end{enumerate}

\subsection{Numerical results}

\noindent We present, in the following, a few numerical results.\\
 \noindent Tests have been made for:

$$ r = 0.06 \quad , \quad \sigma = 0.3 \quad , \quad  E = 100 \quad , \quad  T=60\, \text{days}$$

\noindent The following figures display the graph of the difference $\Delta-\Delta_{classical}$ as a function of stock price $S$ and time $t$, for the values of the shape parameter $\lambda$ such that, at the time $\displaystyle  \frac{T}{2}$, for $S=100$:

\begin{enumerate}
\item[$\rightsquigarrow$] $\Delta-\Delta_{classical}=0.15$, which leads to the choice $\lambda=2.37163$ ;

\item[$\rightsquigarrow$] $\Delta-\Delta_{classical}=0.25$, which leads to the choice $\lambda=1.91079$ ;

\item[$\rightsquigarrow$] $\Delta-\Delta_{classical}=0.35$, which leads to the choice $\lambda=1.67543$ ;

\end{enumerate}

\newpage

\begin{figure}[h!]
 \center{\psfig{height=7cm,width=8.5cm,angle=0,file=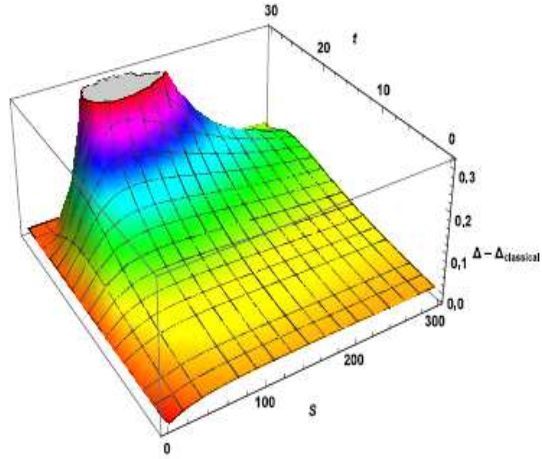}} \\
 \caption{The graph of the difference $\Delta-\Delta_{classical}$ as a function of stock price $S$ and time $t$, for $\lambda=2.37163$.}
  \end{figure}

\begin{figure}[h!]
 \center{\psfig{height=7cm,width=8.5cm,angle=0,file=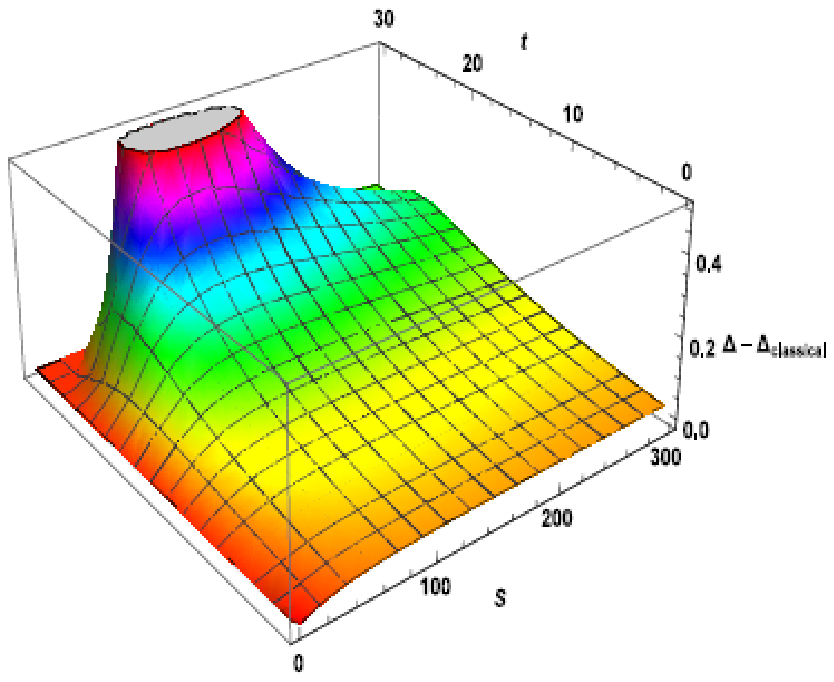}} \\
 \caption{The graph of the difference $\Delta-\Delta_{classical}$ as a function of stock price $S$ and time $t$, for $\lambda=1.91079$.}
  \end{figure}

\begin{figure}[h!]
 \center{\psfig{height=7cm,width=8.5cm,angle=0,file=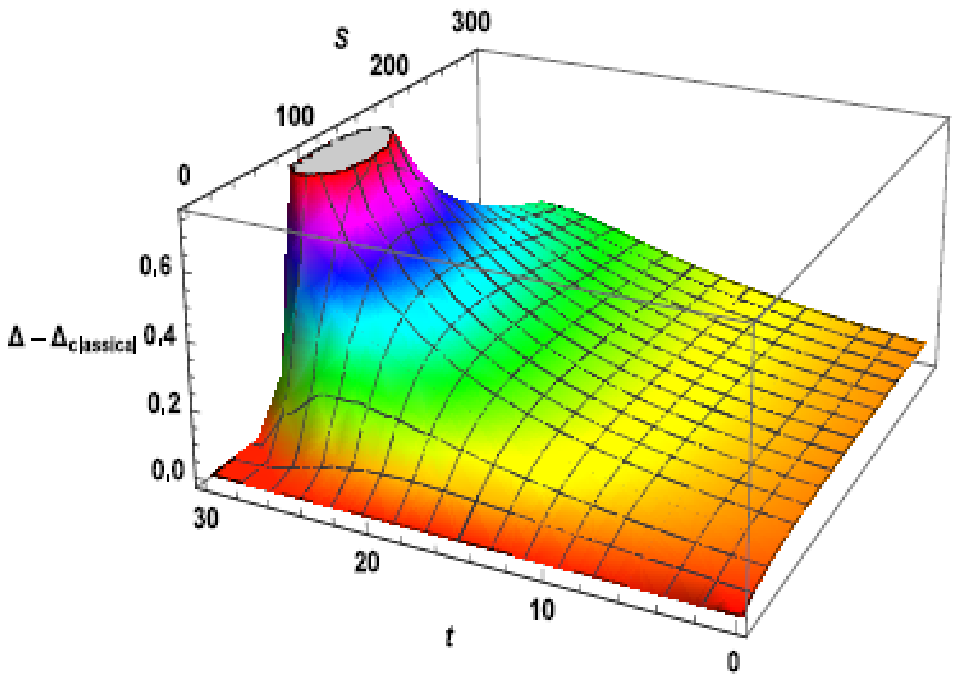}} \\
 \caption{The graph of the difference $\Delta-\Delta_{classical}$ as a function of stock price $S$ and time $t$, for $\lambda=1.67543$.}
  \end{figure}
\vskip 3cm


\end{document}